\documentclass[12pt]{amsart}
\usepackage{epsf, amssymb, xy, latexsym, amsmath, graphics}
\newcommand{\comp}{\, {{\rm o}}\,}
\def \1{^{-1}}
\def \2{^{-2}}
\def \3{^{-3}}
\def \pr{^\prime}
\def \prpr{^{\prime\prime}}
\renewcommand{\phi}{\varphi}

\newcommand{\cf}{{\mathcal F}}

\newcommand{\ch}{{\mathcal H}}

\newcommand{\cl}{{\mathcal L}}
\newcommand{\cm}{{\mathcal M}}

\newcommand{\co}{{\mathcal O}}

\newcommand{\cs}{{\mathcal S}}
\newtheorem{theorem}{Theorem}

\newtheorem{lemma}{Lemma}
\newtheorem{proposition}{Proposition}
\newcommand{\C}{{\mathbb  C}} 
\newcommand{\N}{{\mathbb  N}} 
\newcommand{\Z}{{\mathbb  Z}}

\newcommand{\IP}{{\mathbb P}}
\newcommand{\osum}{\oplus} % I confuse these

\newcommand{\orbfg}{\pi_1^{orb}} %orbifold fundamental group
\newcommand{\ok}{\to}
\newcommand{\eksi}{\backslash}

\newcommand{\onto}{\twoheadrightarrow}
\newcommand{\moins}{\backslash}
\newcommand{\sol}{\bigl\langle}
\newcommand{\sag}{\bigr\rangle}
\newcommand{\supp}{\mbox{supp}}
\newcommand{\ssag}{{\rangle\!\rangle}}
\newcommand{\ssol}{\langle\!\langle}
\newcommand{\lcm}{{\rm lcm}}

\title{Smooth finite abelian uniformizations of projective spaces and Calabi-Yau Orbifolds}
\author{A. Muhammed Uluda{\u g}} 

\begin{document}
\thanks{This work was supported by T{\"U}B{\.I}TAK grant KAR{\.I}YER-104T136}
\dedicatory{Dedicated to Mehmet {\c C}ift{\c c}i}
\begin{abstract}  
We give a  classification of smooth complex manifolds with a finite abelian group action,
such that the quotient is isomorphic to a projective space.  
The case where the manifold is a Calabi-Yau is studied in detail.
\end{abstract} 
\maketitle  
\pagestyle{myheadings}
\markboth{Abelian uniformizations of projective spaces}{A. Muhammed Uluda{\u g}}

\section{Introduction}
Let $M$ be a complex manifold with a faithful action by a finite abelian group $G$ such that $M/G\simeq \IP^n$. 
In this paper we give a complete classification of such pairs $(M,G)$ 
and study the case where $M$ is a Calabi-Yau manifold in detail.
\par 
A well-known example of  a pair $(M, G)$ is a cyclic multiple plane: 
Let 
$$
S:\{[z_0:\dots:z_{n}]\in\IP^{n}\;| \; P(z_0,\dots, z_{n})=0\}
$$
be a smooth hypersurface in $\IP^n$ 
defined by a homogeneous polynomial $P$ of degree $n+2$. 
Then the hypersurface in $\IP^{n+1}$ defined as
$$
M:\{ [z_0:\dots:z_n:z_{n+1}]\in\IP^{n+1}\;|\; z_{n+1}^{n+2}=P(z_0,\dots, z_{n})\} 
$$
is smooth of degree $n+2$, which implies that  $M$ is a Calabi-Yau variety of dimension $n$. 
In dimension two $M$ is a smooth quartic surface 
and in dimension three $M$ is a smooth quintic threefold.  
Let $\omega$ be a primitive $n+2$nd root of unity. 
Then the cyclic group $G:=\Z/(n+2)$ acts on $M$, the action of $i\in G$ is given by  
$$
[z_0:\dots:z_n:z_{n+1}]\in M\to [z_0:\dots:z_n:\omega^iz_{n+1}]\in M
$$
Consider the projection 
$$
\varphi: [z_0:\dots:z_n:z_{n+1}]\in\IP^{n+1} \to [z_0:\dots:z_{n}]\in\IP^n,
$$
its restriction to $M$ is precisely the quotient map $M\ok M/G$, which shows that $M/G$ is $\IP^n$.
Evidently, $G$ fixes the points $M\cap \{ z_{n+1}=0\}$ and the image of this set under $\varphi$ is $H$. 
In other words,
$\varphi:M\ok \IP^n$ is a Galois covering of degree $n+2$, branched along the hypersurface $H$.  
\par
Let  $(M,G)$ be a pair with $M/G\simeq \IP^n$. 
The corresponding projection map $\varphi:M\ok \IP^n$ induces an orbifold
structure $(\IP^n, b_\varphi)$ on $\IP^n$, where  $b_\varphi: \IP^n\ok \N$ is the map sending $p\in \IP^n$ 
to the order of the stabilizer $G_q\subset G$, where $q$ is a point in $\varphi\1(p)$. 
In the example given above, the induced orbifold  is $(\IP^n, b_\varphi)$, where 
$$
b_\varphi(p):=\left\{ 
\begin{array}{ll} 
n+2&p\in S\\
1&p\notin S
\end{array}\right.
$$
The \textit{locus} of an orbifold $(\IP^n,b)$ is defined to be the hypersurface $\supp(b-1)$.
In our case, the locus $(\IP^n,b_\varphi)$ is precisely the hypersurface $S$. 
In Section 2 we give a brief introduction to orbifolds.
The main Theorem~\ref{thm:1} is proved in 2.1.
The orbifold euler number of orbifolds $(\IP^n,b)$ with a linear locus is computed in 2.2.
\par 
An orbifold uniformized by a smooth Calabi-Yau manifold is called a \textit{Calabi-Yau orbifold}.
Section 3 is devoted to the classification of CY orbifolds admitting an abelian uniformization.
In dimension one, this classification is classical. 
K3 orbifolds with a locus of degree $\leq 5$ were classified in \cite{uludagk3}. 
There are no K3 orbifolds with a locus of degree $>6$. 
\par
This work was  partially realized during my stay along September 2005 in I.H.E.S., to which 
I am thankful.

\section{Orbifolds}
Here we give a brief introduction to orbifolds following Kato~\cite{kato}.
For details one may also consult~\cite{mylectures}.
Let $M$ be a connected smooth complex manifold, $G\subset \mbox{Aut}(M)$ a properly
discontinuous subgroup and put $X:=M/G$. 
Then the projection $\phi:M\ok X$ is a branched Galois 
covering  endowing $X$ with a map $b_\phi:X\ok\N$ defined by 
$b_\phi(p):=|G_q|$ where 
$q$ is a point in $\phi^{-1}(p)$ and $G_q$ is the isotropy subgroup of $G$ at $q$. 
The pair $(X,b_\phi)$ is said to be uniformized by $\phi:M\ok (X,b_\phi)$.
An \textit{orbifold} is a pair $(X,b)$ of an irreducible normal analytic space $X$ with a 
function $b:X\ok \N$ such that  the pair $(X,b)$ is locally finitely uniformizable. 
We shall mostly be concerned with the case where $X$ is a smooth complex manifold.
In case $c |b$, the orbifold $(X,c)$ is said to be a \textit{suborbifold} of 
$(X,b)$. Let $\phi:(Y,1)\ok (X,c)$ be a uniformization of $(X,c)$, e.g. $b_\phi=c$. 
Then $\phi:(Y,b\pr)\ok (X,b)$ is called an \textit{orbifold covering}, where 
$b\pr:=b\comp\phi/c\comp\phi$. 
The orbifold $(Y,b\pr)$ is called the 
\textit{lifting of $(X,b)$ to the uniformization $Y$ of $(X,c)$} and  
will be denoted by  $(X,b)/(X,c)$. 

\par
Let $X$ be a smooth manifold. Let $(X,b)$ be an orbifold,  
$H:={\rm supp}(b-1)$ be its \textit{locus}  and let $H_1,\dots,H_r$ be the irreducible components of $H$.
Then $b$ is constant on $H_i\moins {\rm sing}(H)$; so let $m_i$ be this number.
By abuse of language, the divisor $\Sigma_im_iH_i$ is also called the \textit{locus} of $(X,b)$ 
(In fact, the divisor $\Sigma_im_iH_i$ completely determines the $b$-function).  
The \textit{orbifold fundamental group} $\orbfg(X,b)$ of $(X,b)$ is the group defined by
$\orbfg(X,b):=\pi_1(X\moins H)/\ssol\mu_1^{m_1},\dots, \mu_r^{m_r}\ssag$
where $\mu_i$ is a meridian of $H_i$ and $\ssol\ssag$ denotes the normal closure. 
The \textit{local orbifold fundamental group $\orbfg(X,b)_p$ at a point $p\in X$}
is the orbifold fundamental group $\orbfg(\co,p)$ of the restriction of $(X,b)$ to a 
sufficiently small neighborhood of $\co$ of $p$. 
\par
Let $(X,b)$ be an orbifold and $\rho:\orbfg(X,b)\onto G$ be a surjection onto a finite group. 
Then there exists a Galois covering of $M\ok X$ branched along $H$ with $G$ as the Galois group. 
Under which conditions the covering is a uniformization of $(X,b)$ ?
The following theorem (stated in a slightly different terminology in Kato's monograph) answers this question. 

\begin{theorem}\label{uni}
Let $(X,b)$ be an orbifold and $\rho:\orbfg(X,b)\onto G$ a surjection onto a finite group.
The corresponding Galois covering $\phi:M\ok X$ is a uniformization of $(X,b)$ if and only 
if for any $p\in X$, the composition of the maps
\begin{equation}\label{unif}
\orbfg(X,b)_p \stackrel{\iota_p}{\longrightarrow} \pi_1(X, b)\stackrel{\rho}{\longrightarrow} G
\end{equation}
is an injection, where $\iota_p$ is the map induced by the inclusion $(X,b)_p\hookrightarrow (X,b)$. 
\end{theorem}
\par
If the groups $\orbfg(X,b)$ and $G${\nobreakspace}are non-abelian, it is a hopeless task to verify the condition of Theorem~\ref{uni}
However, assuming $G$ to be abelian simplifies matters significantly.

\subsection{Abelian orbifolds on $\IP^n$}
An \textit{abelian orbifold} is an orbifold which admits an abelian uniformization.
Let $(\IP^n,b)$ be an abelian orbifold,  $H:=\supp(b-1)$  its locus, and let $H_1,\dots, H_r$ be the 
irreducible components of $H$. By the following lemma, the $b$-function is very restricted.
\begin{lemma}\label{lem:1}
If $(\IP^n,b)$ is an abelian orbifold, then\\
(i) Any singularity of the locus $H$ is locally of the form 
$\{ z_1z_2\dots z_k=0\, :\, (z_1, \dots z_n) \in \C^n\}$ for some $k\leq n$. \\
(ii) There exists an r-tuple $(m_1,\dots, m_r)\in \N^r$ such that $b$ is given by 
$b(p):=\prod_{p\in H_i} m_i$. \\
(iii) The groups $\pi_1(\IP^n\eksi H)$ and $\orbfg(\IP^n, b)$ are abelian when $n>1$.\\
(iiii) The hypersurfaces $H_i$ are smooth and are in general position. \\

\end{lemma}

\proof
Suppose $(\IP^n,b)$ admits an abelian uniformization and let $G$ be the corresponding abelian Galois group.
By Theorem~\ref{uni}, for any $x\in \IP^n$, the local group $\orbfg(X,b)_x$ injects into $G$. 
Hence all local groups are abelian and all local germs admit abelian uniformizations. 
An abelian covering germ $\C^n_O\ok \C^n_O$ is of the form  
$$
\varphi: (z_1,\dots, z_n)\in \C^n_O \to (z_1^{m_1}, \dots, z_n^{m_n} )\in \C^n_O 
$$
where $m_i\in \Z_{\geq 1}$. 
Hence the singularities of $H$ must be of the form $z_1z_2\dots z_k=0$ for some $k\in[1,n]$. 
This proves (i). Hence locally the $b$-map is of the form $b_\phi(p)=\prod_{p\in L_i}m_i$, where $L_i:=\{z_i=0\}$.
This proves (ii). 
By the Zariski conjecture proved by Fulton-Deligne and by (i) the group 
$\pi_1(\IP^n\eksi H)$ and its quotient $\orbfg(\IP^n, b)$ are abelian, hence (iii) is also proved.
To prove the last claim suppose $p$ is a singular point of $H_1$. 
For simplicity assume that $p$ do not belong to $\cup_2^rH_i$.
Around $p$, the hypersurface $H_1$ is of the form $z_1\dots z_k=0$ for some $1<k\leq n$.
Let $\mu_{1,i}$ be a branch of $z_i=0$ for $i\in [1,k]$. 
Then 
$$
\orbfg(\IP^n,b)_p\simeq \sol \mu_{1,1}, \dots, \mu_{1,k}\; |\; m_1\mu_{1,i}=0 \;\mbox{for} \; i\in k]\sag\simeq  (\Z/(m_1))^k
$$
On the other, since $H_1$ is irreducible the meridians $\mu_{1,i}$ are conjugate elements of the group
$\pi_1(\IP^n\eksi B_b)$. 
Since this latter group is abelian, the image of $\iota_p$ is generated by a single meridian and is cyclic.
Hence $\iota_p$ can not be injective, contradicting Theorem~\ref{uni}. \hfill $\Box$

\bigskip
Lemma~\ref{lem:1} does not give sufficient conditions for the existence of a uniformization. 
Let $(\IP^n, b)$ be an orbifold satisfying the conclusions of Lemma~\ref{lem:1}.
Let $H:=\cup_1^rH_i\subset \IP^n$ be its locus and let 
\begin{equation}\label{eq:1}
b(p):=\prod_{p\in H_i} m_i
\end{equation}
be its $b$-function as required by Lemma \ref{lem:1}, 
where $(m_1,\dots, m_r)\in \N^r$. 
Put 
$$
f_i:=\frac{m_i}{\mbox{gcd}(m_i,d_i)}
$$

\begin{theorem}\label{thm:1}
The orbifold $(\IP^n, b)$ admits a finite abelian smooth uniformization 
if and only if any prime power dividing one among $f_1, \dots, f_r$ divides at least $n$ other $f_i$'s.
\end{theorem}
\proof
The group $\pi_1(\IP^n\eksi H)$ admits the presentation
$$
\pi_1(\IP^n, b)\simeq \sol \mu_1, \dots, \mu_r, \,|\, \sum_1^rd_im_i=0 \sag
$$
where $\mu_i$ is a meridian of $H_i$. 
Hence the group $\orbfg(\IP^n, b)$ admits the presentation
$$
\orbfg(\IP^n, b)\simeq \sol \mu_1, \dots, \mu_r, \,|\, m_1\mu_1=\dots=m_r\mu_r=\sum_1^rd_im_i=0 \sag
$$
This group is of order   
$$
|\orbfg(\IP^n, b)|=\frac{\prod_{i\in[1,r]} m_i}{\lcm\{f_i\, :\, i\in [1,r]\} }
$$
Now let $B\subset [1,r]$ with $|B|\leq n$. 
Let $p$ be a point in $\cap_{i\in B}H_i$, which is not in $\cap_{i\in C}H_i$ for any $C\supseteq B$, in other words,
$$
p\in \bigcap_{i\in B}H_i\eksi \left(\bigcup_{C\supseteq B}\bigcap_{i\in C}H_i\right)
$$
One has the isomorphism of orbifold germs
$$
(\IP^n, b)_p \simeq (\C, 1)^{n-|B|}_0 \times \prod_{i\in [1,n]}(\C, b_i)_0,
$$
where $b_i(0)=m_i$ and $b_i(z)=1$ if $z\neq 0$.
Hence, the local orbifold fundamental group around $p$ is 
$$
\orbfg(\IP^n, b)_p\simeq \sol \{\mu_i\,: i\in B\}\,|\, m_i\mu_i=0 \, : i\in B \sag,
$$
which is of order $\Pi_{i\in B}m_i$. 
Let  $\iota_p:\orbfg(\IP^n, b)_p \ok \orbfg(\IP^n, b)$ be the homomorphism induced by the inclusion. 
One has 
$$
\mbox{Coker}(\iota_p)\simeq\sol \{\mu_i\,: i\notin B\}\,|\, m_i\mu_i=\sum_{i\notin B}d_i\mu_i=0 \sag,
$$
The homomorphism $\iota_p$ is an injection if and only if $\mbox{Coker}(\iota_p)$ is of order
$$
\frac{|\orbfg(\IP^n, b)|}{|\orbfg(\IP^n, b)_p|}=\frac{\prod_{i\notin B}m_i}{\lcm\{f_i\, :\, i\in [1,r]\} }
$$
On the other hand, one has 
$$
|\mbox{Coker}(\iota_p)| = \frac{\prod_{i\notin B} m_i}{\lcm\{f_i\, :\, i\notin B\} }
$$
Therefore, $\iota_p$ is injective if and only if 
\begin{equation}\label{eq:3}
{\lcm\{f_i\, :\, i\in [1,r]\} }={\lcm\{f_i\, :\, i\notin B\} }
\end{equation}
Recall that $B$ is a subset of the interval $[1,r]$ with $|B|\leq n$ elements, and $p\in \cap_{i\in B}H_i$ is a singular
point of the arrangement $\cup_{i\in [1,r]}H_i$ of multiplicity $|B|$. 
The map $\iota_p$ must be injective for any $p\in \IP^n$, so that 
(\ref{eq:3}) must be valid for any $B\subset [1,r]$ with $|B|\leq n$. 
It is easily seen that if  (\ref{eq:3}) is valid 
for any $B$ with  $|B|=n$, then it is valid for any $B$ with $|B|\leq n$. 
But if (\ref{eq:3}) is valid for any $B$ with  $|B|=n$, then the claim of the theorem easily follows \hfill $\Box$

How to interpret the claim of Theorem~\ref{thm:1}? 
Let us first discuss the case $r=1$. 
Thus the locus $H=H_1$ of $(\IP^n,b)$ is a smooth irreducible hypersurface of degree $d_1$, and $f_1=m_1/ \gcd(m_1, d_1)$. 
Now if $f_1\neq 1$, then it has a prime divisor, which must divide at least $n$ other 
$f_i$'s, which is impossible (unless $n=0$). 
Hence $f_1=1$, which implies that $m_1$ is a divisor of $d_1$. 
It is well known that the orbifolds $(\IP^n, b)$ are uniformizable in this case.
\par
Now let us discuss the case $r=2$. 
If $n=1$ then Theorem~\ref{thm:1} implies $f_1=f_2$, which in turn implies that $m_1=m_2$ since $H_i$ are reduced. 
This is the well known orbifold $(\IP^1, b)$, which is uniformized via the power map $\IP^1\ok \IP^1$ of degree $m_1$.
If $n\geq 2$ as in the case of $r=1$ one has $f_1=f_2=1$, in other words $m_i$ divides $d_i$.   
For example, let $H_1$ be a smooth quadric and $H_2$ be a smooth cubic. 
Then $m_1=2$ and $m_2=3$. 
\par
Let $p$ be a prime, $\alpha$ an integer $>0$, and take integers $\alpha_{n+2}, \dots, \alpha_{r} \in [0, \alpha]$. 
Then the following $r$-tuple, 
$$
(p^\alpha, \dots, p^\alpha, p^{\alpha_{n+2}}, \dots, p^{\alpha_{r}})
$$
as well as all its permutations, satisfies the condition of  Theorem~\ref{thm:1}.
Any $r$-tuple $(f_1, \dots, f_r)$ satisfying the condition of  Theorem~\ref{thm:1} admits a unique factorization 
into a product of such $r$-tuples with distinct $p$, under the operation of componentwise multiplication.  
\par
Note that if $H_i$ is a hyperplane then $d_i=1$, so that $f_i=m_i$. 
Hence for hyperplane arrangements the conditions simplifies significantly. 

\subsection{Euler numbers of abelian orbifolds with a linear locus}
Let $(\IP^n, b)$ be an orbifold. 
Then there is an induced stratification $\cs$ of $\IP^n$ with $b$-constant strata.
Let $b(S)$ be the constant value of the $b$-function on the stratum $S\in\cs$.  
The euler number of  $(\IP^n, b)$ is defined as
$$
e_{orb}(\IP^n, b):=\sum_{S \in {\mathcal S}}\frac{e(S)}{b(S)}
$$
so that, if $\phi:M\ok (\IP^n,b)$ is a uniformization with $G$ as the Galois group, then 
$e(M)=|G|e_{orb}(\IP^n, b)$. 
\par
Now let $(\IP^n, b)$ be an orbifold whose locus is a hyperplane arrangement $\ch:=\{H_1, \dots, H_r\}$. 
Let $\cl(\ch)$ be its intersection lattice. 
Then the induced stratification is 
$$
\cs=\{   {\cm(\ch^A )} \;: \; A \in \cl \},
$$
where as usual $ \ch^A$ denotes the restriction of $\ch$ to $A$ and 
$\cm(\ch)$ denotes the complement. 
Let $(m_1, \dots, m_r)$ be as in Lemma~\ref{lem:1}. 
Then $b(H_i)=m_i$. 
The constant value assumed by the $b$-function along an element $A\in \cl$ is the product $\prod_{A\subset H_i}b(H_i)$. 
Hence one has 
$$
e_{orb}(\IP^n, b)=\sum_{A\in \cl } \frac{e(\cm(\ch^A))}{\prod_{A\subset H_i}m_i}
$$
The number $e(\cm(\ch^A))$ can be computed as follows. 
Suppose $A$ is of rank $k$ in $\cl$ (i.e. $A$ is of codimension $k$ in $\IP^n$).
Then $A$ lies in the intersection of exactly $k$ hyperplanes $H_i$.   
By the genericity assumption, $\ch^A$ is an arrangement of $r-k$ hyperplanes 
in general position in  $A\simeq \IP^{r-k}$. 
Let $e(r-k,n-k)$ be the euler number of $\cm(\ch^A)$.
Let $\ch$ be an arrangement of hyperplanes in general position and $H\in \ch$. 
As usual, denote by $\ch_H$ the deletion $\ch\eksi \{H\}$. 
The equation
$$
\cm(\ch)=\cm(\ch_H) -  \cm(\ch^{H})
$$
yields the recursion 
$$
e(r,n)=e(r-1, n)-e(r-1, n-1)
$$
Tabulation of $e(r,n)$ gives the following Pascal triangle with alternating signs.

\medskip
$$
\begin{array}{r|rrrrrrrrr}
  &0&1&2& 3& 4 &  5&  6&    7&   8\\
\hline
\IP^0&1&1&1& 1& 1 &  1&  1&    1&   1\\
\IP^1&2&1&0&-1&-2& -3& -4&  -5&  -6\\
\IP^2&3&1&0& 0& 1 &  3&  6&  10& 15\\
\IP^3&4&1&0& 0& 0 & -1& -4&-10&-20\\
\IP^4&5&1&0& 0& 0 &  0&  1&    5& 15 \\
\IP^5&6&1&0& 0& 0 &  0&  0&   -1&  -6\\
\IP^6&7&1&0& 0& 0 &  0&  0&    0&   1\\
\IP^7&8&1&0& 0& 0 &  0&  0&    0&   0
\end{array}
$$

\medskip
Hence $e(r,n)=(-1)^n\left(\begin{array}{c} r-2\\ n \end{array}\right)$, 
in other words  
$$
e(\cm(\ch^A))=e(r-k,n-k)=(-1)^{n-k}\left(\begin{array}{c} r-2-k\\ n-k \end{array}\right) 
$$
This gives
$$
e_{orb}(\IP^n, b)=\sum_{k=0}^n (-1)^{n-k}\left(\begin{array}{c} r-2-k\\ n-k \end{array}\right)
\sum_{\scriptsize{\begin{array}{c} B\subset [1,r]\\ |B|=k\end{array}}}\prod_{i\in B}\frac{1}{m_i}
$$
The change of parameters $s_i:=1-1/m_i$ gives, after recollecting terms the final formula 
in the following proposition.

\begin{proposition}
The orbifold euler number of $(\IP^n,b)$ is
\begin{equation}
e_{orb}(\IP^n, b)=\sum_{j=0}^n (-1)^j (n+1-j) 
\sum_{\scriptsize{\begin{array}{c} B\subset [1,r]\\ |B|=j\end{array}}}\prod_{i\in B} s_i
\end{equation}
\end{proposition}

\section{Calabi-Yau orbifolds}
Let $M$ be a uniformization of the orbifold $(\IP^n,b)$ and let $\phi:M\ok \IP^n$ be the corresponding 
covering map.  Denote by $D_\phi$ the ramification divisor of $\phi$ and let $\sum m_iH_i$ be the locus of $(\IP^n,b)$.
One has 
$$
\begin{array}{ccl}
K_M&=&\phi^*K_{\IP^n}+D_\phi\\
&=&\phi^*K_{\IP^n}+\sum_i\frac{m_i-1}{m_i} \phi^*H_i\\
&=&\phi^*\Bigl(K_{\IP^n}+\sum_i\bigl(1-\frac{1}{m_i}\bigr)H_i\Bigr)
\end{array}
$$
Now one has $K_m\sim -(n+1)L$ and $H_i \sim d_iL$, where $L$ is the class of a hyperplane in $\IP^n$. 
Hence $K_M$ is trivial if and only if the equality
$$
\sum_i s_i=\sum_i\bigl(1-\frac{d_i}{m_i}\bigr)=n+1
$$
is satisfied.
As the following trivial lemma shows, this is a very restrictive condition on $(\IP^n,b)$. 
\begin{lemma}\label{lem:2}
If $(\IP^n,b)$ is a Calabi-Yau orbifold  with a locus $H:=\cup_1^rH_i$ of degree $d=\sum d_i$, then $n+2\leq d\leq 2n+2$.
Moreover, if $d=2n+2$, then $m_1=m_2=\cdots =m_{r}=2$.
\end{lemma}

Hence in any dimension there are finitely many families of abelian Calabi-Yau orbifolds on $\IP^n$, 
which can be effectively classified. For $n=1$ this classification is classical.
We give the details of this classification as a guide to the classifications in dimension 2 and 3, 
tabulated in the coming pages.

\bigskip
$$ 
\begin{array}{l|llrrlll}
\IP^1&d&\mbox{orbifold}&e&|\orbfg|&\delta & \mbox{sub-orbifolds and coverings}\\
\hline\hline
{\bf 1}&2&[\infty, \infty]&0&\infty&0&[m,m]{\bf 1}\\
{\bf 2}&3&[2,2,\infty]&0&\infty&0&[2,2,1]{\bf 1}, [1,2,2]{\bf 2} \\
{\bf 3}&3&[2,3,6]&0&\infty&0&[2,1,2]{\bf 5}, [1,3,3]{\bf 6}\\
{\bf 4}&3&[2,4,4]&0&\infty&0&[1,2,2]{\bf 6}, [1,4,4]{\bf 6}, [2,2,1]{\bf 4}\\ 
{\bf 5}&3&[3,3,3]&0&\infty&0&[3,3,1]{\bf 5}\\
{\bf 6}&4&[2,2,2,2]&0&\infty&1&[2,2,1]{\bf 6}\\
\hline
\end{array} 
$$ 

\bigskip
\noindent
Here are some remarks on how to understand the tables:

\noindent
\begin{itemize}
\item In dimension $n$, the orbifold $[\infty,\infty  \dots, \infty ]$ ($n+1$ times) is uniformized by $\C^n$, via a 
multi-exponential map.
Strictly speaking, this is not a CY orbifold. 
It is included in the tables in order to complete the picture.
\item \underline{$d$}: The degree of the locus of the orbifold
\item \underline{orbifold}: In dimension 1, the notation 
$[m,m\pr ,\dots ]$ means that the locus consists of distinct points $p, p\prpr,\dots, $, 
such that $b(p)=m$, $b(p\pr)=m\pr$, etc. 
In higher dimensions, the notation $[m_{d}, m\pr_{d\pr}, \dots]$ means  that the locus consists of 
hypersurfaces $H$ of degree $d$, $H\pr$ of degree $d\pr$, etc, 
and the generic value the $b$-function takes on $H$ is $m$, on $H\pr$ is $m\pr$, etc. 
In case the corresponding hypersurface is linear, $m_1$ is denoted simply by $m$.
\item \underline{$e$}: The euler number of the universal uniformization (which is simply connected). 
Thus $e=e_{orb}|\orbfg|$.
\item \underline{$|\orbfg|$}: 
The order of the orbifold fundamental group. 
Hence this is the degree of the universal uniformization.
\item \underline{$\delta$}: The dimension of the family.
\item \underline{sub-orbifolds and coverings}: 
In dimension $n$, the linear orbifold $[m, m, \dots, m]$ ($n+1$ times) is uniformized by $\IP^n$, 
for any $m$. 
Some of the CY orbifolds admits these as sub-orbifolds and 
consequently can be lifted to their uniformization. 
The lifting is another abelian CY orbifold on $\IP^n$.  
The bold number adjacent to a sub-orbifold in the table is an internal reference to the line 
containing the lifting of the orbifold to the uniformization of this sub-orbifold. 
For example, in dimension 1, the orbifold $[1,3,3]$ is a sub-orbifold of $[2,3,6]$, and the lifting of 
$[2,3,6]$ to the uniformization of $[1,3,3]$ is an orbifold of type $[2,2,2,2]$. 
Another example in dimension 2 is the 3-dimensional family $[2_2,3,3,3]$, 
which admit $[1,3,3,3]$ as a sub-orbifold. 
The lifting of this family to the uniformization of $[1,3,3,3]$ is a 3-dimensinal 
sub-family of the 19-dimensional family $[2_6]$, whose loci are defined by equations $P[x^3:y^3:z^3]=0$, 
the polynomial $P$ being homogeneous of degree 2. 
\item The euler numbers of  non-linear orbifolds which appear as coverings of linear orbifolds can be computed 
by using the fact that their universal uniformizations are isotopic so that their euler numbers are same.
\item In dimensions $\geq 4$, the tables includes only CY orbifolds with linear loci.
This is for sake of brevity only, there are many non-linear cases.
\end{itemize}

\subsection{A non-abelian CY orbifold}
Some abelian CY orbifolds have non-trivial automorphism groups. 
Their quotients sometimes yields CY orbifolds with a non-abelian uniformizing group.
In dimension 2, many examples were constructed in~\cite{uludagk3}.
In dimension 3 consider the CY orbifold $[2,2,2,2,2,2,2,2]$. 
Suppose that its locus consists of the hyperplanes $z_0\pm z_1\pm z_2\pm z_3=0$. 
This orbifold is invariant under the action of the group $\Z/(2)\osum\Z/(2)\osum \Z/(2)$
such that the element$(i,j,k)$ acts by 
$$
[z_0,z_1,z_2,z_3]\longrightarrow [z_0, (-1)^iz_1, (-1)^jz_2, (-1)^kz_3]
$$
The quotient of $\IP^3$ under this action is again $\IP^3$, via the map
$$
\phi_2([z_0,z_1,z_2,z_3]):=[z_0^2,z_1^2,z_2^2,z_3^2]
$$
The quotient of  $[2,2,2,2,2,2,2,2]$  under the same group action is therefore the non-abelian CY orbifold 
$(\IP^3, b)$, whose locus is given by the equation 
$$
\{z_0z_1z_2z_3(\sqrt{z_0}+\sqrt{z_1}+\sqrt{z_2}+\sqrt{z_3})=0\} \subset \IP^3
$$
The $b$-function takes the generic value 2 on the locus. 
Computing the degree of the universal uniformization of  $(\IP^3, b)$ we find that 
the group $\orbfg(\IP^3, b)$ is a non-abelian group of order $2^7\times2^3=1024$.
This latter orbifold is invartiant under the symmetric group on four letters $\{z_0, z_1, z_2, z_3\}$.
However, this time the quotient of $\IP^3$ is a singular space. 
Anyway, this shows that the uniformizing CY threefold has a 
group of automorphisms of order $1024\times 4!=24576$.

\subsection{CY manifolds of Enriques type}
An abelian orbifold may have several uniformizations, in other words 
the universal uniformization is not necessarily the only uniformization.
There are also some CY orbifolds admitting intermediate uniformizations.
For any $n$, the orbifold fundamental group of the CY orbifold $[2,2,\dots, 2]$ ($2n+2$ times) 
admits the presentation 
$$
\orbfg([2,2,\dots, 2])\simeq\sol \mu_1, \dots, \mu_{2n+2} \, |\, 2\mu_1=\dots=2\mu_{2n+2}=\Sigma_1^{2n+2}\mu_i=0\sag
$$
Let $\langle \alpha \rangle $ be the subgroup generated by $\alpha:=\mu_1+ \dots+ \mu_{n+1}$ and let $G$ be the quotient by this subgroup.
It has the presentation
$$
G\simeq    \sol \mu_1, \dots, \mu_{2n+2} \, |\, 2\mu_1=\dots=2\mu_{2n+2}=\Sigma_1^{2n+2}\mu_i=\Sigma_1^{n+1}\mu_i=0\sag
$$
This group satisfies the conditions of Theorem~\ref{uni}-therefore the corresponding 
uniformization exists. 
Since $\alpha $ is of order $2$, this is not the universal uniformization, and its fundamental group is $\Z/(2)$. 
In dimension 2 this uniformization is an Enriques surface. 
The element $\alpha$ can be chosen in 
$\frac{1}{2}{\scriptsize \left(\begin{array}{c}2n+2\\n+1\end{array}\right)}$ different ways, 
hence there are as many intermediate uniformizations. 

\subsection{CY orbifolds and configuration spaces}
Families of K3 surfaces frequently has a nice interpretations as quotients of symmetric spaces
(see~\cite{dolgachev}).
Perhaps the simplest example is the following: 
Consider the family $\cf$ of K3 orbifolds of type $[6_2, 3,3]$,  
whose locus consists of a conic $Q$ with two lines in general position. 
Now consider $Q$ as a projective line $\IP^1$ and the intersection points of $Q$ with the lines 
as a coloured configuration of points on $\IP^1$, such that the points lying on the same line has the same color.
Hence $\cf$ can be interpreted as a configuration space of coloured points on $\IP^1$, 
and is isogeneous to the configuration space of 4 points on $\IP^1$. 
The moduli space of 4-tuples on $\IP^1$ is a quotient of the Poincar{\'e} disc in the well-known, classical manner.
Note that tangent lines to $Q$ lifts to the uniformization of a K3 orbifold of type $[6_2, 3,3]$
as a pair of elliptic fibrations. 
In the examples we discuss below this feature is pertinent, except the last one.
The fibrations corresponding to the orbifolds $[6_2, 3,3]$ and $[2_2, 3,3,3]$ are isotrivial.
\par
Consider now the family $\cf$ of K3 orbifolds of type $[2_2, 3,3,3]$ (or of type $[4_2,2,2,2]$),  
whose locus consists of a conic $Q$ with three lines in general position.
This family can be interpreted as a configuration space of coloured points on $\IP^1$, 
and is isogeneous to the configuration space of 6 points on $\IP^1$. 
The moduli space of 6-tuples on $\IP^1$ is known to be a ball-quotient (see~\cite{demo}). 
Similarly, the family of K3 orbifolds of type $[2_2,2,2,2,2]$ is related to the moduli space of 
8 points on $\IP^1$, which is a ball quotient. It was shown in~\cite{yoshida} that the 
family of K3 orbifolds of type $[2,2,2,2,2,2]$ is a quotient of a type IV-symmetric domain.
It is of interest to know whether other families of Calabi-Yau orbifolds constructed in this paper 
admit similar interpretations.
{}
\bigskip
\tiny

\noindent
\begin{center} Abelian K3 orbifolds \end{center}
$$ 
\begin{array}{l|llrrlll}
\IP^2&d&\mbox{orbifold}&e&|\orbfg|&\delta & \mbox{sub-orbifolds and coverings}\\
\hline\hline
{\bf 1}&3&[\infty,\infty,\infty]&0&\infty&0&[m,m,m]{\bf 1}\\
{\bf 2}&4&[2,6,6,6]&24&72&0&[2,2,2,1]{\bf 4}, [1,2,2,2]{\bf 6}, [1,3,3,3]{\bf 12}, [1,6,6,6]{\bf 14}\\
{\bf 3}&4&[4,4,4,4]&24&64&0&[2,2,2,1]{\bf 7}, [4,4,4,1]{\bf 5} \\
{\bf 4}&4&[6_2,3,3]&24&18&1&\\
{\bf 5}&4&[4_4]&24&4&6&\\
{\bf 6}&5&[2_2, 3,3,3]&24&18&3&[1,3,3,3]{\bf 14}\\
{\bf 7}&5&[4_2,2,2,2]&24&16&3&[1,2,2,2]{\bf 5}\\
{\bf 8}&5&[3_3, 2_2]&24&6&6&\\
{\bf 9}&6&[2,2,2,2,2,2]&24&32&4&[2,2,2,1,1,1]{\bf 10}\\
{\bf 10}&6&[2_2,2,2,2,2]&24&32&5&[1,1,2,2,2]{\bf 13}\\
{\bf 11}&6&[2_2, 2_2, 2_2]&24&8&7&\\
{\bf 12}&6&[2_3, 2,2,2]&24&8&7&[1,2,2,2]{\bf 14}\\
{\bf 13}&6&[2_4, 2_2]&24&4&11&\\
{\bf 14}&6&[2_6]&24&2&19&\\
\hline
 \end{array}
 $$

\newpage
 \noindent
\begin{center} Abelian Calabi-Yau 3-orbifolds \end{center}
 $$
 \begin{array}{l|llrrlll}
\IP^3&\!\! d&\!\! \mbox{orbifold}&\!\! e&\!\! |\orbfg|&\!\! \delta &\!\!  \mbox{sub-orbifolds and coverings}\\
\hline\hline
{\bf 1}&\!\! 4&\!\! [\infty,\! \infty,\! \infty,\! \infty]&\!\! 0&\!\! \infty&\!\! 0&\!\! [m,\! m,\! m,\! m]{\bf 1}\\
{\bf 2}&\!\! 5&\!\! [2,\! 5,\! 10,\! 10,\! 10]&\!\! -288&\!\! 1000&\!\! 0&\!\! 
[1,\! 5,\! 5,\! 5,\! 5]{\bf 30},[2,\! 1,\! 2,\! 2,\! 2]{\bf 6}\\
{\bf 3}&\!\! 5&\!\! [2,\! 8,\! 8,\! 8,\! 8]&\!\! -296&\!\! 1024&\!\! 0&\!\! 
[1,\! 2,\! 2,\! 2,\! 2]{\bf 12},[1,\! 4,\! 4,\! 4,\! 4]{\bf 27},[1,\! 8,\! 8,\! 8,\! 8]{\bf 34},[2,\! 2,\! 2,\! 2,\! 1]{\bf 7}\\
{\bf 4}&\!\! 5&\!\! [3,\! 6,\! 6,\! 6,\! 6]&\!\! -204&\!\! 648&\!\! 0&\!\! 
[1,\! 2,\! 2,\! 2,\! 2]{\bf 12},[1,\! 3,\! 3,\! 3,\! 3]{\bf 22}, [1,\! 6,\! 6,\! 6,\! 6]{\bf 19}, [3,\! 3,\! 3,\! 3,\! 1]{\bf 16}\\
{\bf 5}&\!\! 5&\!\! [5,\! 5,\! 5,\! 5,\! 5]&\!\! -200&\!\! 625&\!\! 0&\!\! 
[1,\! 5,\! 5,\! 5,\! 5]{\bf 8}\\
{\bf 6}&\!\! 5&\!\! [5_2,\! 5,\! 5,\! 5]&\!\! -200&\!\! 125&\!\! 3\\
{\bf 7}&\!\! 5&\!\! [8_2,\! 4,\! 4,\! 4]&\!\! -200&\!\! 128&\!\! 3\\
{\bf 8}&\!\! 5&\!\! [5_5]&\!\! -200&\!\! 5&\!\! 55\\
{\bf 9}&\!\! 6&\!\! [2,\! 2,\! 3,\! 3,\! 6,\! 6]&\!\! -120&\!\! 216&\!\! 3&\!\! 
[1,\! 1,\! 3,\! 3,\! 3,\! 3]{\bf 28},[2,\! 2,\! 1,\! 1,\! 2,\! 2]{\bf 14}\\
{\bf 10}&\!\! 6&\!\! [2,\! 2,\! 4,\! 4,\! 4,\! 4]&\!\! -176&\!\! 256 &\!\! 3&\!\! 
[1,\! 1,\! 2,\! 2,\! 2,\! 2]{\bf 26},[1,\! 1,\! 4,\! 4,\! 4,\! 4]{\bf 32}, [2,\! 2,\! 2,\! 2,\! 1,\! 1]{\bf 15}\\
{\bf 11}&\!\! 6&\!\! [3,\! 3,\! 3,\! 3,\! 3,\! 3]&\!\! -144&\!\! 243&\!\! 3&\!\! 
[1,\! 1,\! 3,\! 3,\! 3,\! 3]{\bf 18}\\
{\bf 12}&\!\! 6&\!\! [2_2,\! 4,\! 4,\! 4,\! 4]&\!\! -296&\!\! 128&\!\! 6&\!\! [1,\! 4,\! 4,\! 4,\! 4]{\bf 34}, [1,\! 2,\! 2,\! 2,\! 2]{\bf 27}\\
{\bf 13}&\!\! 6&\!\! [3_2,\! 3,\! 3,\! 3,\! 3]&\!\! -204&\!\! 81&\!\! 6&\!\!  [1,\! 3,\! 3,\! 3,\! 3]{\bf 19}\\
{\bf 14}&\!\! 6&\!\! [3_2,\! 3_2,\! 3,\! 3]&\!\! -120&\!\! 27&\!\! 9\\
{\bf 15}&\!\! 6&\!\! [4_2,\! 4_2,\! 2,\! 2]&\!\! -176&\!\! 32&\!\! 9\\
{\bf 16}&\!\! 6&\!\! [6_3,\! 2,\! 2,\! 2]&\!\! -204&\!\! 24&\!\! 13\\
{\bf 17}&\!\! 6&\!\! [2_2,\! 4_4]&\!\! &\!\! 8&\!\! 28 \\
{\bf 18}&\!\! 6&\!\! [3_3,\! 3_3]&\!\! -144&\!\! 9&\!\! 23\\
{\bf 19}&\!\! 6&\!\! [3_6]&\!\! -204&\!\! 3&\!\! 83\\
{\bf 20}&\!\! 7&\!\! [4_2,\! 2,\! 2,\! 2,\! 2,\! 2]&\!\! &\!\! 64&\!\! 9 &\!\!  [1,\! 2,\! 2,\! 2,\! 2]{\bf 17}\\
{\bf 21}&\!\! 7&\!\! [2_2,\! 4_2,\! 2,\! 2,\! 2]&\!\! -176&\!\! 32&\!\! 12\\
{\bf 22}&\!\! 7&\!\! [3_3,\! 2,\! 2,\! 2,\! 2]&\!\! -204&\!\! 24&\!\! 16&\!\!   [1,\! 2,\! 2,\! 2,\! 2]{\bf 19}\\
{\bf 23}&\!\! 7&\!\! [2_2,\! 2_2,\! 3_3]&\!\! &\!\! 12&\!\! 22\\
{\bf 24}&\!\! 7&\!\! [3_3,\! 2_4]&\!\! &\!\! 38 &\!\! 6 \\
{\bf 25}&\!\! 8&\!\! [2,\! 2,\! 2,\! 2,\! 2,\! 2,\! 2,\! 2]&\!\! -128&\!\! 128&\!\! 9&\!\!   [1,\! 1,\! 1,\! 1,\! 2,\! 2,\! 2,\! 2]{\bf 28}\\
{\bf 26}&\!\! 8&\!\! [2_2,\! 2_2,\! 2,\! 2,\! 2,\! 2]&\!\! -176&\!\! 32&\!\! 15&\!\!   [1,\! 1,\! 2,\! 2,\! 2,\! 2]{\bf 32}\\
{\bf 27}&\!\! 8&\!\! [2_4,\! 2,\! 2,\! 2,\! 2]&\!\! -296&\!\! 16&\!\! 31&\!\!   [1,\! 2,\! 2,\! 2,\! 2]{\bf 34}\\
{\bf 28}&\!\! 8&\!\! [2_2,\! 2_2,\! 2_2,\! 2_2]&\!\! -128&\!\! 16&\!\! 21\\
{\bf 29}&\!\! 8&\!\! [2_3,\! 2_3,\! 2,\! 2]&\!\!-120&\!\! 8&\!\! 29\\
{\bf 30}&\!\! 8&\!\! [2_5,\! 2,\! 2,\! 2]&\!\! -288&\!\! 8&\!\! 49\\
{\bf 31}&\!\! 8&\!\! [2_4,\! 2_2,\! 2_2]&\!\! &\!\! 8&\!\! 37\\
{\bf 32}&\!\! 8&\!\! [2_4,\! 2_4]&\!\! -176&\!\! 4&\!\! 53\\
{\bf 33}&\!\! 8&\!\! [2_6,\! 2_2]&\!\! &\!\! 4&\!\! 76\\
{\bf 34}&\!\! 8&\!\! [2_8]&\!\! -296&\!\! 2&\!\! 164\\
\hline
\end{array}
$$

\newpage
\begin{center} Some higher dimensional linear abelian Calabi-Yau orbifolds \end{center}
\tiny
$$
\begin{array}{l|llrrl}
&degree&orbifold&&\\
\hline\hline
\IP^4&&&&\\ 
\hline
&d=6&[2,10,10,10,10,10]&&\\
&&[6,6,6,6,6,6]&&\\
&d=7&[2,2,3,6,6,6,6]&&\\
&&[2,4,4,4,4,4,4]&&\\
&d=10&[2,2,2,2,2,2,2,2,2,2]&&\\
\hline\hline
\IP^5&&&&\\ 
\hline
&d=7&[2,7,14,14,14,14,14]&&\\
&&[2,12,12,12,12,12,12]&&\\
&&[3,4,12,12,12,12,12]&&\\
&&[3,9,9,9,9,9,9]&&\\
&&[4,8,8,8,8,8,8]&&\\
&&[7,7,7,7,7,7,7]&&\\
&d=8&[2,2,6,6,6,6,6,6]&&\\
&&[2,3,3,6,6,6,6,6]&&\\
&&[4,4,4,4,4,4,4,4]&&\\
&d=9&[2,2,2,3,3,3,6,6,6]&&\\
&&[2,2,2,4,4,4,4,4,4]&&\\
&&[3,3,3,3,3,3,3,3,3]&&\\
&d=12&[2,2,2,2,2,2,2,2,2,2,2,2]&&\\
\hline\hline
\IP^6&&&&\\ 
\hline
&d=8&[2,14,14,14,14,14,14,14]&&\\
&&[8,8,8,8,8,8,8,8]&&\\
&d=9&[2,3,6,6,6,6,6,6,6]&&\\
&d=10&[2,2,2,3,3,3,6,6,6,6,6]]&&\\
&d=14&[2,2,2,2,2,2,2,2,2,2,2,2,2,2]&&\\
\hline\hline
\IP^7&&&&\\ 
\hline
&d=9&[2,9,18,18,18,18,18,18,18]&&\\
&&[2,16,16,16,16,16,16,16,16]&&\\
&&[3,5,15,15,15,15,15,15,15]&&\\
&&[3,12,12,12,12,12,12,12,12]&&\\
&&[4,6,12,12,12,12,12,12,12]&&\\
&&[5,10,10,10,10,10,10,10,10]&&\\
&&[9,9,9,9,9,9,9,9,9]&&\\
&..& .... &&\\
\end{array}
$$

\end{document}